\documentclass{article}
\usepackage{amsmath,amssymb,amsfonts,epsfig}

\newcommand{\Gt}{\overset{\sim}{\mathcal{G}}}
\newcommand{\Got}{\overset{\sim}{G_1}}
\newcommand{\Gtt}{\overset{\sim}{G_2}}
\newcommand{\FUM}{\overset{\sim}{FUM}}

\title{On some projective unitary qutrit gates}
\author{Claire Levaillant}
\begin{document}
\maketitle

\textbf{Abstract.}  As part of a protocol, we braid in a certain way six anyons of topological charges $222211$ in the Kauffman-Jones version of $SU(2)$ Chern-Simons theory at level $4$. The gate we obtain is a braid for the usual qutrit $2222$ but with respect to a different basis. With respect to that basis, the Freedman group of \cite{LEV} is identical to the $D$-group $D(18,1,1;2,1,1)$. We give a physical interpretation for each Blichfeld generator of the group $D(18,1,1;2,1,1)$. Inspired by these new techniques for the qutrit, we are able to make new ancillas, namely $\frac{1}{\sqrt{2}}(|1>\,+|3>)$ and $\frac{1}{\sqrt{2}}(|1>\,-|3>)$, for the qubit $1221$. %Furthermore, by braiding and measurement we are able to make a new gate which is not a braid for the qubit $1221$.
\section{Setting}

$\;\;$ Recently, there has been some interest in finding qutrit gates which are universal for quantum computation. When the group of qutrit gates in the projective unitaries $PU(3)$ acts irreducibly on $\mathbb{C}^3$, a result of \cite{F}
provides a sufficient condition named by the authors condition $(\star)$ for an $SU(3)$-subgroup of single projective unitary qutrit gates to form a dense set of $PU(3)$. This condition finds its origins in a $2002$ work \cite{FKL} by Michael Freedman, Alexei Kitaev and Jacob Lurie. An older result from Jean-Luc Brylinski and Ranee Brylinski \cite{BB} implies that such a dense set of $1$-qutrit gates together with a $2$-qutrit entangling gate is universal for quantum computation. Therefore, there have been some attempts and hopes, starting from a finite group of projective unitary qutrit gates obtained by anyonic braiding, to add an extra projective unitary gate which would this time be obtained by braiding and interferometric measurement and would make the group become infinite. We believe that such a group would then satisfy to the conditions mentioned above for density.\\
\indent In \cite{BL}, we study a finite subgroup of $SU(3)$ arising from anyonic braiding. This group has order $162$ and is later enlarged to a group of order $648$, the Freedman group, by a fusion operation (FFO for future reference) due to Mike Freedman, see \cite{LEV}. Both groups, the one of order $162$ and its extension of order $648$ are isomorphic to $D$-groups in the $1916$ classification of finite $SU(3)$-subgroups by Blichfeld (later augmented with two new groups), namely to $D(9,1,1;2,1,1)$ and to $D(18,1,1;2,1,1)$ respectively. In \cite{LEV}, it is shown further that the $D$-group $D(18,1,1;2,1,1)$ is the Freedman group, with respect to a different basis, that is both groups are conjugate.
Classically and originally, the group $D(18,1,1;2,1,1)$ is defined by three matrix generators which first appeared in the $1916$ book by Blichfeld as part of the three generic generators for the groups $D(n,a,b;d,r,s)$ from the series $(D)$. Our paper introduces a new set of four generators for the $D$-group $D(18,1,1;2,1,1)$, but the group is only generated by three of them. These generators all arise from anyonic braiding and FFO. %By adding a fourth generator, we are able to easily retrieve the fact that the Freedman group is isomorphic to $D(18,1,1;2,1,1)$ without the need of a computer program like GAP or without exhibiting an isomorphism between both groups like in \cite{LEV}. Moreover, we are able to easily enlighten a physical interpretation for the original three generators of $D(18,1,1;2,1,1)$.
Both the Freedman group of order $648$ and our physical interpretation of $D(18,1,1;2,1,1)$ contain the center $Z_3$ of $SU(3)$, hence we note that the number of projective unitary qutrit gates available to us remains the same.
In this first part of the paper, we consider the qutrit $2222$ and a pair of $1$'s, do some specific braids and fail to obtain a new gate. Of course the number of protocols available to us is extremely large, so our failure does not imply that by choosing such an ancilla we won't ever obtain an interesting gate by braiding and measurement. Two fundamental facts are enlightened from this first part. First, when doing a full twist $\sigma_2$ on four particles $2211$, it results in swapping the topological charges $0$ and $2$. Second, when doing a single braid $\sigma_2$ on four particles $2211$, we obtain a qubit $2121$ with the same proportion of $|1>$ and $|3>$. Since doing $\sigma_1$ braids only introduces phases, we can thus make a qubit $1221$ with equal norms of $|1>$ and $|3>$. This was unknown fact in \cite{F} where in some protocols using braiding and interferometric measurement on the qubit $1221$, we were missing such ancillas which play a crucial role for the no-leakage condition. \\
%Last, by considering next the qubit $1221$ and this time a pair of $2$'s as ancilla, we succeed to make an additional gate for the qubit $1221$ in a way that is described in $\S\,4$.

% but that this new way of making qutrit ancilla could happen to be very useful for making anyonic gates beyond braiding.
%But we have a new way of obtaining projective unitary qutrit gates by braiding and FFO.
%We think that some of these gates could possibly be used in a fruitful manner for producing ancillas which could be used for finding $PU(2)$ qubit gates obtained by braiding and forced measurement.
\section{Result}
We state below our result.
\newtheorem{Theorem}{Theorem}
\begin{Theorem}
The group $\Gt$ generated by the four matrices
$$\begin{array}{l}\begin{array}{cc}\overset{\sim}{G_1}=\begin{pmatrix}
&&e^{\frac{7i\pi}{9}}\\&-e^{\frac{4i\pi}{9}}&\\e^{\frac{7i\pi}{9}}&&\end{pmatrix}
&\overset{\sim}{G_2}=\begin{pmatrix}&e^{\frac{7i\pi}{9}}&\\e^{\frac{7i\pi}{9}}&&\\&&-e^{\frac{4i\pi}{9}}
\end{pmatrix}\end{array}\\
\begin{array}{cc}
\overset{\large\sim}{FUM}=\begin{pmatrix}
-e^{\frac{2i\pi}{3}}&&\\&-e^{\frac{2i\pi}{3}}&\\&&e^{\frac{2i\pi}{3}}
\end{pmatrix}
&N=\begin{pmatrix} -e^{-\frac{i\pi}{9}}&&\\&-e^{-\frac{i\pi}{9}}&\\&&e^{i\frac{2\pi}{9}}\end{pmatrix}
\end{array}
\end{array}$$
is a finite subgroup of $SU(3)$ of order $648$. It is isomorphic to a semi-direct product $C_{6}\times C_{18}\rtimes S_3$.
The generators above are up to phase obtained by the following unitary operations in the Kauffman-Jones version of $SU(2)$ Chern-Simons theory at level $4$.

\begin{center}
\epsfig{file=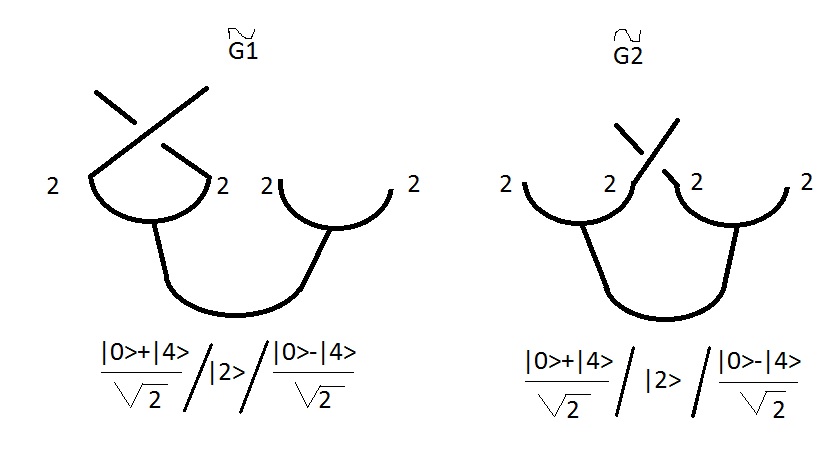, height=5cm}
\end{center}
\begin{center}
\epsfig{file=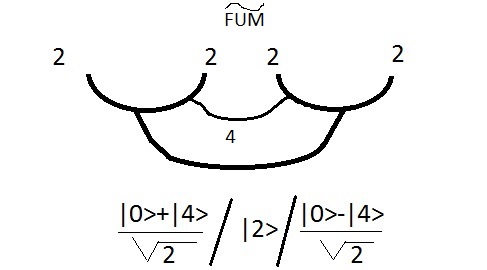, height=4cm}
\end{center}
\begin{center}
\epsfig{file=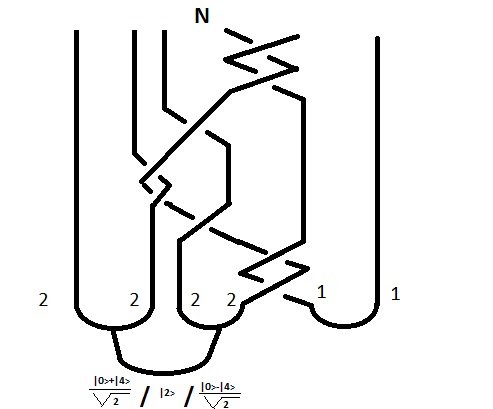, height=9cm}
\end{center}
The generator $N$ belongs to the subgroup generated by $\Gtt$. Moreover, we have
$$\Gt=<\Got,\Gtt,\FUM>=D(18,1,1;2,1,1)$$
\end{Theorem}
\section{Protocol}
A starting point are braids on four anyons of topological charge $2$ in the Jones-Kauffman version of $SU(2)$ Chern-Simons theory at level $4$. We recall below the matrices $G_2$ for a $\sigma_2$-braid and $G_1$ for a $\sigma_1$-braid, also commonly called $R$-matrix, taken from \cite{BL}. All the matrices are defined in $SU(3)$, that is they are defined up to phase.
$$\begin{array}{cc}G_1=\begin{pmatrix}e^{\frac{7i\pi}{9}}&&\\&-e^{\frac{4i\pi}{9}}&\\
&&-e^{\frac{7i\pi}{9}}\end{pmatrix},
&G_2=\begin{pmatrix}
-\frac{1}{2}e^{\frac{4i\pi}{9}}&\frac{e^{\frac{7i\pi}{9}}}{\sqrt{2}}&\frac{1}{2}
e^{\frac{4i\pi}{9}}\\&&\\
\frac{e^{\frac{7i\pi}{9}}}{\sqrt{2}}&0&\frac{e^{\frac{7i\pi}{9}}}{\sqrt{2}}\\
&&\\\frac{1}{2}e^{\frac{4i\pi}{9}}&\frac{e^{\frac{7i\pi}{9}}}{\sqrt{2}}&
-\frac{1}{2}e^{\frac{4i\pi}{9}}
\end{pmatrix}\end{array}$$

\begin{center}
\epsfig{file=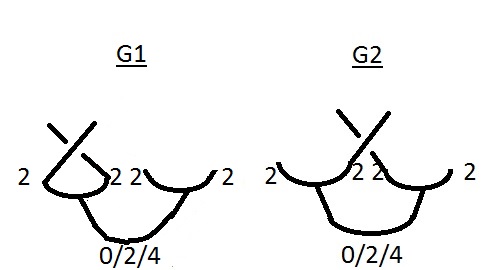, height=4cm}
\end{center}

On the matrices above, we notice the special roles played by the qutrits $|0>$ and $|4>$ on the one hand and $|2>$ on the other hand. Explicitly, braiding anyons $1$ and $2$ maps the qutrit $|0>$ to itself and the qutrit $|4>$ to the qutrit $-|4>$, up to a common phase. Notice further that
$$G_2(|2>)=e^{\frac{7i\pi}{9}}\;\;\frac{|0>+|4>}{\sqrt{2}}$$
and $$G_2\bigg(\frac{|0>+|4>}{\sqrt{2}}\bigg)=e^{\frac{7i\pi}{9}}\;|2>$$
and $$G_2\bigg(\frac{|0>-|4>}{\sqrt{2}}\bigg)=-e^{\frac{4i\pi}{9}}\;\frac{|0>-|4>}{\sqrt{2}}$$
From now on, we will work in the new basis $(e_1,e_2,e_3)$ with
$$e_1=\frac{|0>+|4>}{\sqrt{2}},\;\;e_2=|2>,\;\;e_3=\frac{|0>-|4>}{\sqrt{2}}$$
The matrices of the $\sigma_1$ and $\sigma_2$ braids with respect to this new basis are the following. Again, in all what follows, we write the matrices involved with determinant $1$, that is we drop a phase. And so we get:
$$\begin{array}{cc}\overset{\sim}{G_1}=\begin{pmatrix}
&&e^{\frac{7i\pi}{9}}\\&-e^{\frac{4i\pi}{9}}&\\e^{\frac{7i\pi}{9}}&&\end{pmatrix}
&\overset{\sim}{G_2}=\begin{pmatrix}&e^{\frac{7i\pi}{9}}&\\e^{\frac{7i\pi}{9}}&&\\&&-e^{\frac{4i\pi}{9}}
\end{pmatrix}\end{array}$$
With respect to our new basis, the FFO whose effect is to swap the qutrits $|0>$ and $|4>$
\begin{center}\epsfig{file=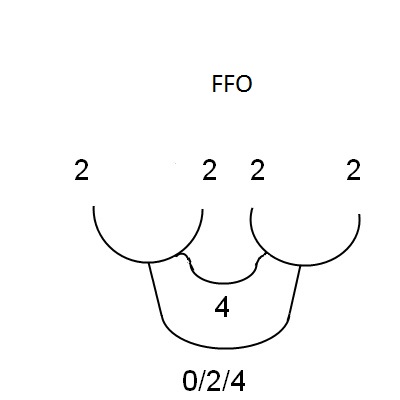, height=6cm}\end{center}
is encoded as follows.
$$\overset{\large\sim}{FUM}=\begin{pmatrix}
-e^{\frac{2i\pi}{3}}&&\\&-e^{\frac{2i\pi}{3}}&\\&&e^{\frac{2i\pi}{3}}
\end{pmatrix}$$
Note in $PU(3)=SU(3)/Z_3$, this matrix is simply
$$\begin{pmatrix}
-1&&\\&-1&\\&&1
\end{pmatrix}$$
We have an analogue for the qubit without need of fusing any particles but simply by using braids. Namely, a full twist like on the figure below has the effect of swapping the qubits $|0>$ and $|2>$. This is a fundamental observation in the protocol we will soon describe.

\begin{center}
\epsfig{file=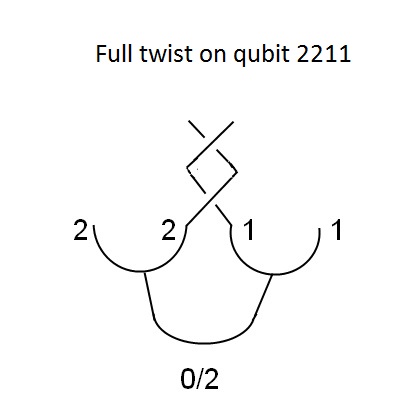, height=6cm}
\end{center}

\textsc{Proof.} The fact that $|0>$ is mapped to $|2>$ essentially relies on the following two points.
\begin{itemize}
\item The quantum dimensions of particles of topological charge $1$ and $3$ are the same.\\
\item The two diagonal coefficients of the squared $R$-matrix $R(2,1)$ are opposite.
\end{itemize}
It then follows that $|2>$ is mapped to $|0>$ by unitarity of the matrix. \\
Let us justify the first point in more details.
Acting on the qubit $|0>$, after doing an $F$-move with horizontal charge line $0$ at the level of the second and the third anyon, followed by two $R$-moves, we obtain the diagram

\begin{center}
\epsfig{file=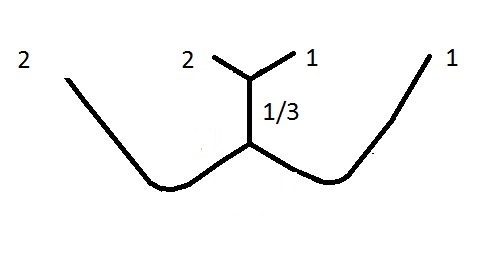, height=4cm}
\end{center}

We then do an $F$-move again. When looking for the $|0>$ projection, the two unitary $6j$-symbols which are involved each contain a $"0"$ which makes them be unitary theta symbols. Using the notations of \cite{ZW} and \cite{LEV}, the two values $\theta^u(1,2,1)$ and $\theta^u(1,2,3)$ are identical since the quantum dimensions of particles of respective topological charge $1$ and $3$ are the same.\hfill $\square$ \\
\begin{center}\underline{\textit{Protocol}}\\\end{center}
$$\begin{array}{l}\\\end{array}$$

\noindent $(i)$ Take a qutrit $2222$ and a pair of $1$'s out of the vacuum. \\Number the anyons $1$, $2$, $3$, $4$, $5$, $6$, those from the qutrit being numbered first. \\\\
$(ii)$ Prepare the qutrit in one of the states $|2>$ or $\frac{|0>+|4>}{\sqrt{2}}$ or $\frac{|0>-|4>}{\sqrt{2}}$.\\\\
$(iii)$ Do a full twist on anyons $4$ and $5$ to "create" a $2$ charge line in between the qutrit and the pair of ones.\\\\
$(iv)$ Use this "extended" version of the qutrit to make braids in such a way that the outcome is a qutrit $2222$ on the first $3$ anyons and a qubit $2211$ on the last $3$.\\\\
$(v)$ Go back to the original configuration of qutrit $2222$ and pair of $1$'s by doing a full twist between anyons $4$ and $5$.\\\\

\noindent Step $(ii)$ is summarized in the following figure and Step $(iii)$ is represented in the figure below it. \\\\
\epsfig{file=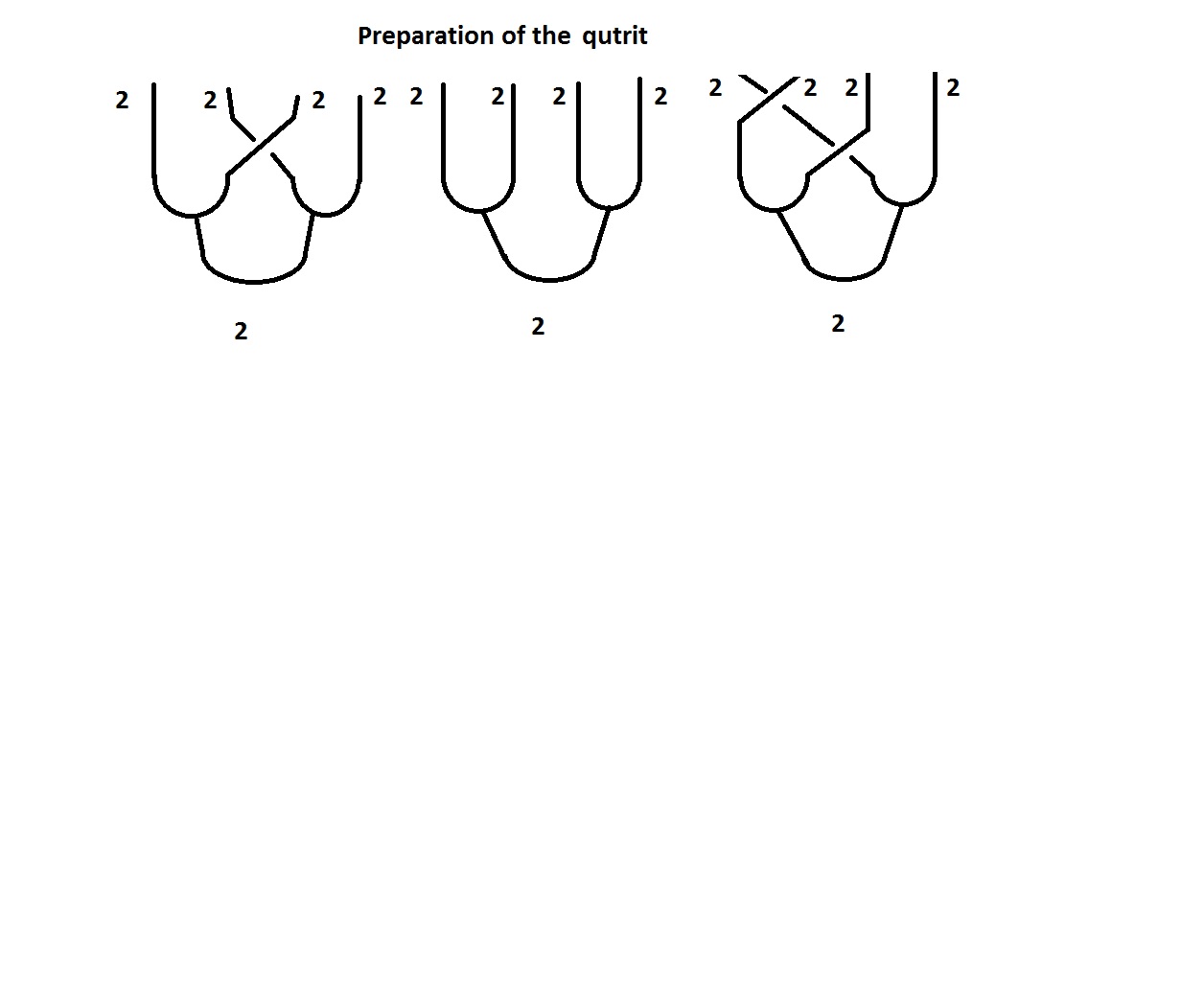, height=10cm}
\vspace{-6cm}
\begin{center}
$\;\;\;\;$\epsfig{file=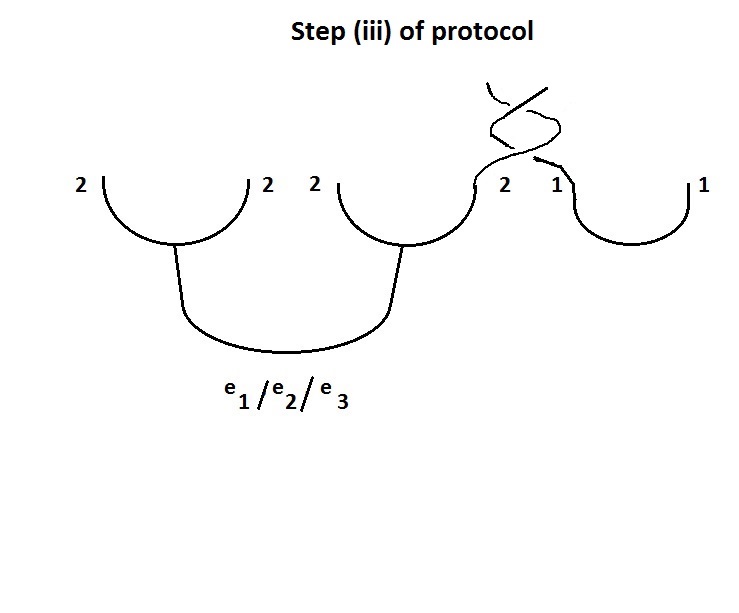, height=6cm}
\end{center}

The braids of step $(iv)$ are now described below.

\begin{center}
\epsfig{file=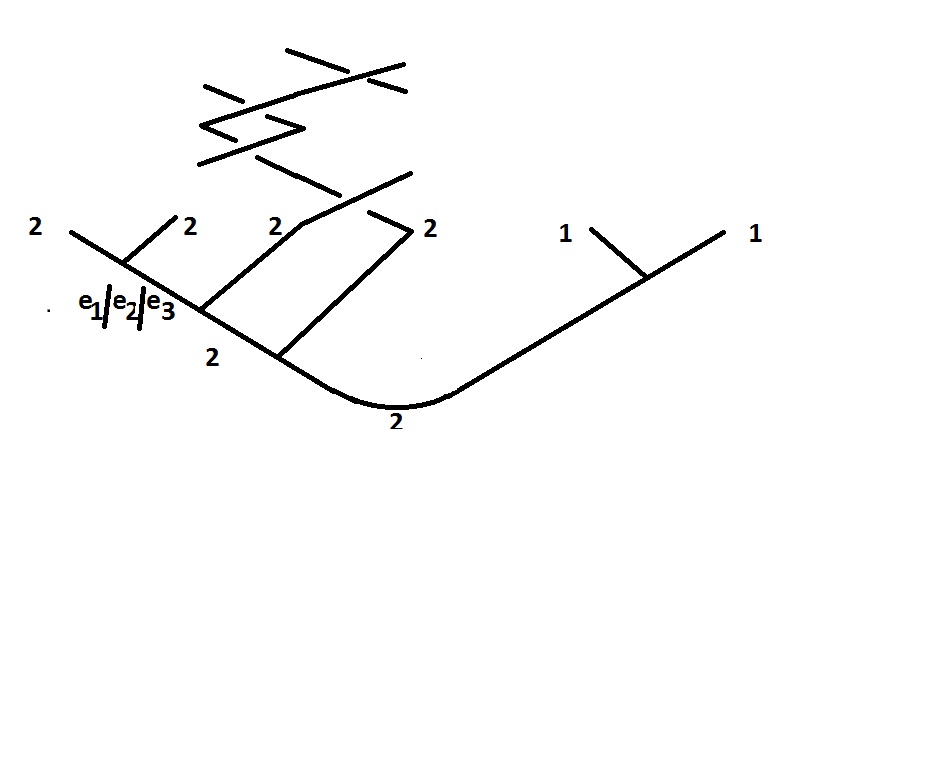, height=10cm}
\end{center}
\vspace{-4cm}
It is a consequence of the fusion rules that a "middle braid" on particles $4222$ or $0222$ or their respective vertical mirror images will map $\mathbb{C}|2>$ into $\mathbb{C}|2>$. \\Moreover, the braiding simply introduces the same phase $e^{\frac{4i\pi}{3}}$, whether dealing with $4222$ or with $0222$.

\begin{center}
\epsfig{file=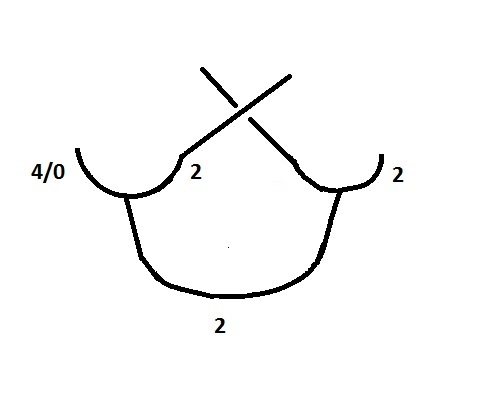, height=3.5cm}\epsfig{file=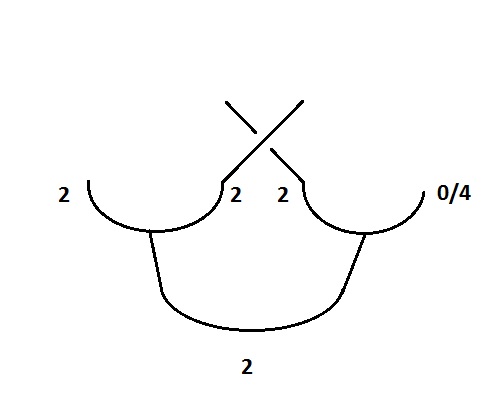, height=3.5cm}
\end{center}

In light of this, it makes sense to do a full twist between anyons $2$ and $3$ on the figure above. It namely allows the charge line adjacent to the input (when going up the tree towards the root) to carry the charge 2 at the end of the braiding process in order for step $(v)$ to be successful independently from the input. After completing the whole protocol, we obtain a new matrix in $SU(3)$, namely
$$N=\begin{pmatrix}
e^{\frac{8i\pi}{9}}&&\\
&e^{\frac{8i\pi}{9}}&\\
&&e^{\frac{2i\pi}{9}}
\end{pmatrix}$$

\section{Group structure}
We will show the following result.
\begin{Theorem}
The group $\overset{\sim}{\mathcal{G}}$ generated by the matrices $\overset{\sim}{G_1}$, $\overset{\sim}{G_2}$, $\overset{\sim}{FUM}$ and $N$ has order $648$ and is isomorphic to a semi-direct product $C_{6}\times C_{18}\rtimes S_3$ with respect to conjugation, for the action provided in Lemma $1$ below. Moreover, it is the group $D(18,1,1;2,1,1)$.
\end{Theorem}
\textsc{Proof.} There is in $\Gt$ a normal subgroup, say $\Delta$, generated by all the diagonal matrices. Moreover, there is a Klein group inside $\Delta$ generated by the two matrices $(\overset{\sim}{FUM})^3$ and its $G_1$-conjugate $\overset{\sim}{G_1}(\overset{\sim}{FUM})^3\overset{\sim}{G_1}^{-1}$.
Indeed, we have
$$\begin{array}{cc}(\FUM)^3=\begin{pmatrix}
-1&&\\&-1&\\&&1
\end{pmatrix},&\Got\,(\FUM)^3\,\Got^{-1}=\begin{pmatrix}
1&&\\&-1&\\&&-1
\end{pmatrix}\end{array}$$

\newtheorem{Lemma}{Lemma}
\begin{Lemma}\hfill\\\\
Our group $\Gt=\bigg< \Got,\,\Gtt,\,\FUM,\,N\bigg>$ is isomorphic to
\hspace{-2cm}
$$\big(<N^2\Got^2>\times\,<N>\,\times\,<(\FUM)^3,\Got(\FUM)^3\Got^{-1}>\big)\;\rtimes\,S_3$$
with $$S_3=\bigg\lbrace \Got^9,\,\Gtt^9,\,\Got^9\Gtt^9,\,\Gtt^9\Got^9,\,\Gtt\Got\Gtt\bigg\rbrace$$
Denoting the latter set by $\lbrace t_3,\,t_1,\,c_1,\,c_2,\,t_2\rbrace$ and the generators from the direct product of two cyclic groups $C_{6}\times C_{18}$ by
$$\begin{array}{cc}\left|\begin{array}{ccc}
x_6&=&N^2\Got^2\,(\FUM)^3\\
x_{18}&=&N\,\Got(\FUM)^3\Got^{-1}
\end{array}\right.&,\end{array}$$

\noindent a presentation for this group is given by
$$\Bigg <x_6,x_{18},t_1,t_2| \begin{array}{l} t_1^2=t_2^2=1=x_6^6=x_{18}^{18}=[x_6,x_{18}]=(t_1t_2)^3,\\
\begin{array}{cc}\begin{array}{ccc}t_1\,x_6\,t_1&=&x_6^{-1},\\
t_1\,x_{18}\,t_1&=&x_6^3\;x_{18}\;\,,\end{array}&\begin{array}{ccc}
t_2\,x_6\,t_2&=&x_6^5\;x_{18}^3\\
t_2\,x_{18}\,t_2&=&x_6^4\;x_{18}^{13}\end{array}\end{array}\end{array}
\Bigg>$$
In the semi-direct product above, $N^2\Got^2$ is the matrix $B^2$ of \cite{BL} with respect to the basis $\bigg(e_1=\frac{|0>+|4>}{\sqrt{2}},\,e_2=|2>,\,e_3=\frac{|0>-|4>}{\sqrt{2}}\bigg)$.\\\\
\noindent The GAP ID for the presentation given above is
$$[648,259]$$
That is our group is the $259$-th group of order $648$ in the SmallGroups library by H. Besche, B. Eick and E. O'Brien dating from the beginning of the $2000$ millenium. This is the same GAP ID as the one of $D(18,1,1;2,1,1)$. \\The group $\Gt$ is precisely the group $$D(18,1,1;2,1,1)$$ \noindent defined by matrix generators by Blichfeld in $1916$.
\end{Lemma}
\textsc{Proof of Lemma.} We look for more cyclic groups generated by diagonal matrices and whose mutual intersections and intersection with the Klein group are trivial.\\ Begin obviously with the subgroup of $\Gt$ generated by the matrix $N$. Notice also $\Got$ squared is a diagonal matrix. We have
$$\begin{array}{cc}\Got^2=\begin{pmatrix}e^{-\frac{4i\pi}{9}}&&\\&e^{\frac{8i\pi}{9}}&\\&&e^{-\frac{4i\pi}{9}},\end{pmatrix}&
N=\begin{pmatrix} e^{\frac{8i\pi}{9}}&&\\&e^{\frac{8i\pi}{9}}&\\&&e^{\frac{2i\pi}{9}}\end{pmatrix}\end{array}$$
Now stare at these matrices. Both matrices have order $9$. Because their diagonal phases in position $(2,2)$ are identical, we see that the two subgroups $<\Got^2>$ and $<N>$ intersect non-trivially only for the $k$-th powers of the generators with $k$ satisfying to $1\leq k\leq 8$ and
$$2k\,\equiv\, -4k\,\equiv\,8k\;\;(mod\,18)$$
This implies that $3$ must divide $k$. Then $k=3$ or $k=6$. In order to solve this unpleasant issue, we must "mix" the generators instead.
We have $$N^2\Got^2=\begin{pmatrix}
e^{\frac{4i\pi}{3}}&&\\&e^{\frac{2i\pi}{3}}&\\&&1\end{pmatrix}$$
and $$\bigg(N^2\Got^2\bigg)^2=\begin{pmatrix}
e^{\frac{2i\pi}{3}}&&\\&e^{\frac{4i\pi}{3}}&\\&&1
\end{pmatrix}$$
And so, we have $$<N^2\Got^2>\;\cap\;<N>=\lbrace I_3\rbrace$$

%We deduce from this discussion that $$<\Got^4>\cap\,<N>=\lbrace I_3\rbrace$$
\noindent In the Klein group, all the elements have order $2$ and in a cyclic group of odd order, all the elements have an odd order. Hence
$<N^2\Got^2>$ and $<N>$ don't intersect with the Klein group. \\

We now exhibit a symmetric group $S_3$ inside $\Gt$. It suffices to notice that
$\Got^9\,\Gtt^9$ and $\Gtt^9\,\Got^9$ are the two usual permutation matrices associated with the respective two cycles of $Sym(3)$. On the other hand, we have
$$\begin{array}{cc}\Got^9=\begin{pmatrix}
&&-1\\&-1&\\-1&&\end{pmatrix},&\Gtt^{9}=\begin{pmatrix}
&-1&\\
-1&&\\
&&-1
\end{pmatrix},\end{array}$$ $$\Gtt\Got\Gtt=\begin{pmatrix}
-1&&\\&&-1\\&-1&
\end{pmatrix}$$

\noindent These matrices provide the additional matrices respectively associated with the three transpositions $(13)$, $(12)$ and $(23)$ of $Sym(3)$.

\noindent It remains to show that each of the $\Gt$-generators $N$, $\overset{\sim}{FUM}$, $\Got$ and $\Gtt$ can be written as a product of an element of the direct product and a group element of $S_3$. The result from the Lemma will then classically follow. \\
First and foremost, we are able to write, using the fact that $\Got$ has order $18$,
$$\Got=N^{-10}\,\bigg(N^2\Got^2\bigg)^5\,\Got^9$$

%since $\Got$ has order $9$, we can write
%$$\Gtt=(\Got^4)^2\,(\Got\Gtt)$$

\noindent Next, it suffices to notice that
$$(\FUM)^4=N^3$$
and so, $$\FUM=N^3(\FUM)^{-3}$$
In particular, we see that the matrix corresponding to the FFO is in the direct product. This was expected since it is a diagonal matrix.
\noindent Further, we have
$$N^2\Gtt^2=(\FUM)^2$$
\noindent We derive
$$\Gtt^2=N^{-2}(\FUM)^2$$

\noindent Now write
$$\Gtt=\bigg(\Gtt^2\bigg)^5\Gtt^9$$
in order to conclude.\\
Finally, it is straightforward to see that $\Gt=D(18,1,1;2,1,1)$. Recall below the Blichfeld generators of $D(18,1,1;2,1,1)$.
$$\begin{array}{l}
F(18,1,1)=\begin{pmatrix} e^{\frac{i\pi}{9}}&&\\&e^{\frac{i\pi}{9}}&\\&&e^{\frac{-2i\pi}{9}}\end{pmatrix}\\\\
\begin{array}{cc}E=\begin{pmatrix} 0&1&0\\0&0&1\\1&0&0\end{pmatrix}&\overset{\sim}{B}=\begin{pmatrix} -1&0&0\\0&0&-1\\0&-1&0\end{pmatrix}
\end{array}\end{array}$$
We see that $E=\Gtt^9\Got^9$ and $\overset{\sim}{B}=\Gtt\,\Got\,\Gtt$, hence $E$ and $\overset{\sim}{B}$ both belong to $\Gt$. Further, we have
\begin{equation}
F(18,1,1)=(\FUM)^{3}\,N^{-1}
\end{equation}
Thus, we see that $F(18,1,1)$ also belongs to $\Gt$ and the Blichfeld generator $F(18,1,1)$ can be expressed in terms of the FFO matrix and the $N$ gate. We conclude that the groups $\Gt$ and $D(18,1,1;2,1,1,)$ are identical since by \cite{LEV} and the current work, they have the same order. We now state below a theorem about a physical interpretation of the original Blichfeld generators of $D(18,1,1;2,1,1)$.
\begin{Theorem}
The Blichfeld generators from $D(18,1,1;2,1,1)$ can be physically realized as follows.

\begin{center}
\epsfig{file=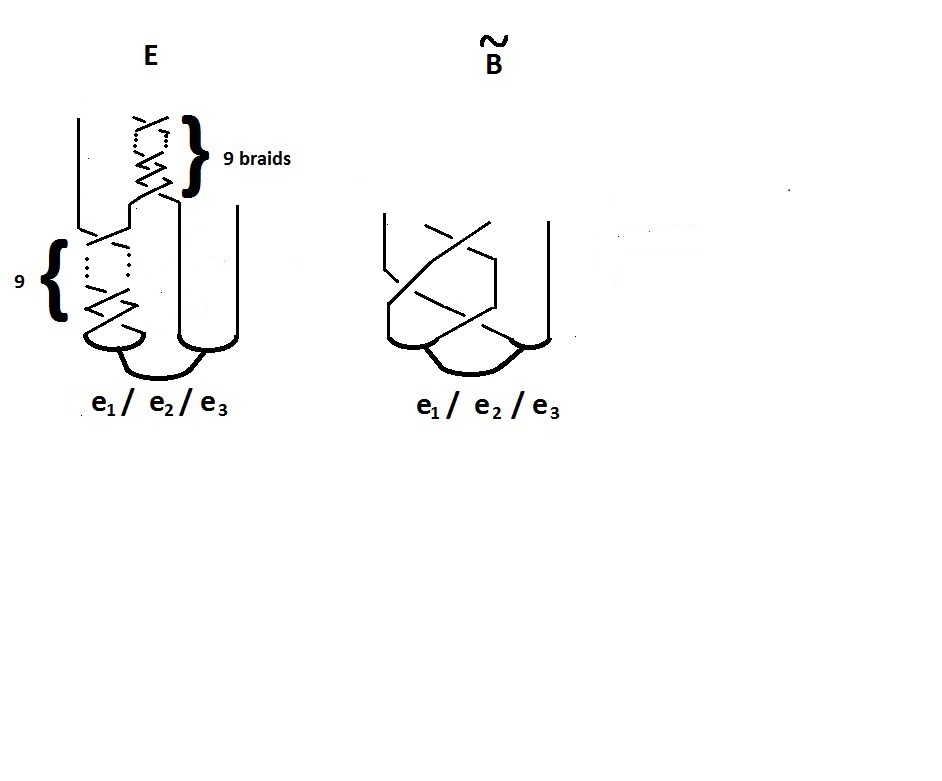, height=9cm}
\end{center}

\vspace{-4cm}

\begin{center}
\epsfig{file=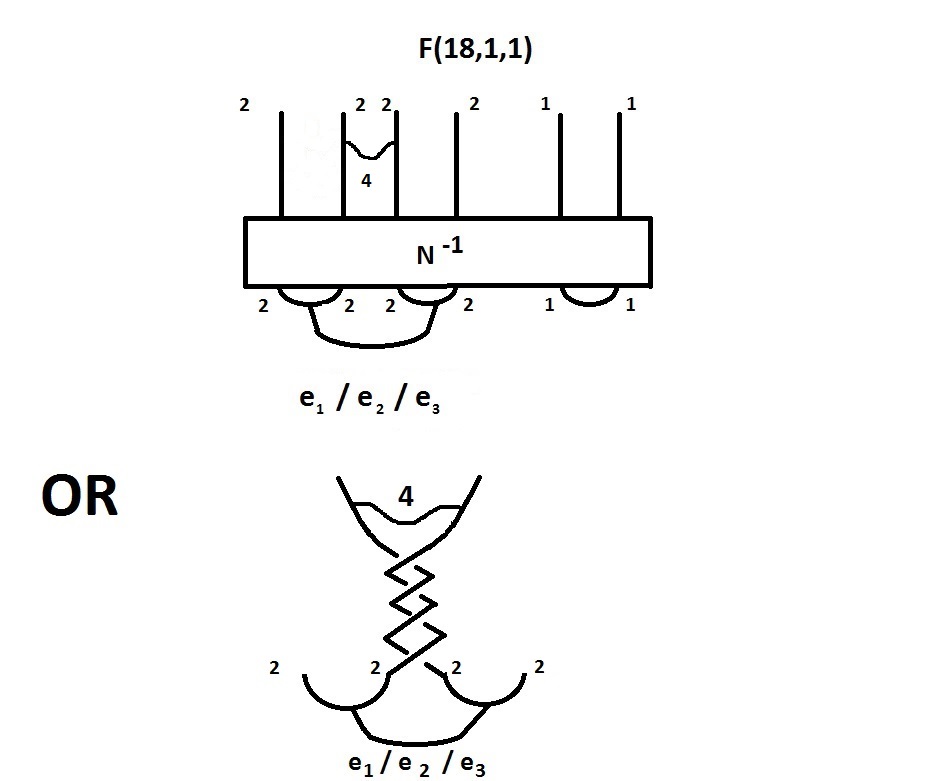, height=8cm}
\end{center}

\end{Theorem}
\noindent The subgroup $<\Got,\Gtt,\FUM>$ has order $648$ since it is conjuage to the Freedman group $<G_1,G_2,FUM>$ of the same order. Hence it is actually the whole group $\Gt$ since as part of our work we showed that $\Gt$ has order $648$. Then the matrix $N$ must be obtained by braiding and FFO. In fact, it is simply obtained by braiding as shown on the figure above. We summarize our results in the Theorem below.
\begin{Theorem}\hfill\\\\
$(i)$ $\qquad\qquad\qquad\qquad\qquad\qquad\Gt\,=\,<\Got,\Gtt,\FUM>$\\\\
$(ii)$ The matrix $N$ is obtained by braiding in an adequate way $4$ anyons $2222$ with respect to the basis $(e_1,e_2,e_3)$. Explicitly, we have
$$N=\Gtt^{-4}$$
\end{Theorem}
\textsc{Proof.} Point $(i)$ was already discussed. As for point $(ii)$, simply notice that
$$N=F(9,1,1)^4$$
and $$F(9,1,1)^{-1}=\Gtt\,t_1$$
Recall $$t_1=\Gtt^9$$ Hence,
$$F(9,1,1)=\Gtt^8$$\hfill $\square$\\

\noindent Last, we comment on the two groups $Fr(162\times 4)$ and $\Gt=D(18,1,1;2,1,1)$.
By \cite{LEV},
$$O^{T}\,Fr(162\times 4)O=\Gt$$
with $$O=\begin{pmatrix} 1/\sqrt{2}&0&-1/\sqrt{2}\\0&1&0\\1/\sqrt{2}&0&1/\sqrt{2}\end{pmatrix}\\$$
where $Fr(162\times 4)$ denotes the Freedman group. We read that $O$ is the transition matrix from
$$(|0>,|2>,|4>)$$ to $$\bigg(\frac{|4>+|0>}{\sqrt{2}},|2>,\frac{|4>-|0>}{\sqrt{2}}\bigg)$$
Thus, we see that $\Gt=D(18,1,1;2,1,1)$ encodes the $\sigma_1$ and $\sigma_2$- braids and FFO on $4$ anyons of topological charge $2$, with respect to either basis
$$\begin{array}{l}
(\frac{|4>+|0>}{\sqrt{2}},|2>,\frac{|4>-|0>}{\sqrt{2}})\\\\
(\frac{|0>+|4>}{\sqrt{2}},|2>,\frac{|0>-|4>}{\sqrt{2}})
\end{array}$$
%Since the matrices of $D(18,1,1;2,1,1)$ are all of the form a permutation matrix times a diagonal matrix with phases on the diagonal, multiplying the last row and last column by $-1$ does not change the matrix. Hence the matrices with respect to either basis are identical.
%Denoting by $\overset{\sim}{O}$ the transition matrix from $(|0>,|2>,|4>)$ to the basis which we called $(e_1,e_2,e_3)$, we obtain
%$$Fr(162\times 4)=\overset{\sim}{O}O^{-1}\,Fr(162\times 4)\,O\overset{\sim}{O}$$
%after noting that $\overset{\sim}{O}^{-1}=\overset{\sim}{O}$.
%From there, since $$\overset{\sim}{O}\,O^{-1}=O\,\overset{\sim}{O}=\begin{pmatrix} 0&0&1\\0&1&0\\1&0&0\end{pmatrix},$$
%we derive that
In other words, we have
$$Fr(162\times 4)=\begin{pmatrix} 0&0&1\\0&1&0\\1&0&0\end{pmatrix}Fr(162\times 4)\begin{pmatrix} 0&0&1\\0&1&0\\1&0&0\end{pmatrix}$$
That is, if we swap the first row and third row and the first column and the third column of a Freedman matrix, we again obtain a Freedman matrix.

\section{New ancilla for the qubit $1221$}

In \cite{F}, we seek ancillas of the form $x\,|1> + y\,|3>$ with $|x|=|y|$ for the qubit $1221$.

\begin{center}\epsfig{file=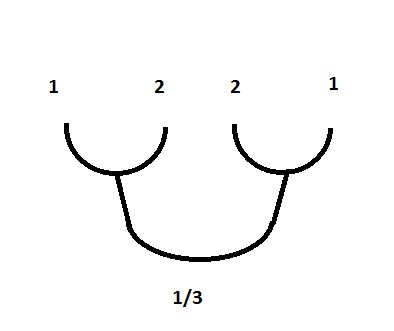, height=4cm}\end{center}
The fact that the norms in $|1>$ and $|3>$ are equal is a necessary condition for no-leakage on some protocols we test which use a combination of braiding and interferometric measurements. Such an ancilla cannot be realized by a combination of $\sigma_1$- and $\sigma_2$-braids on the qubit $1221$. Indeed, the matrix for a $\sigma_2$-braid is the following.
$$\begin{pmatrix}
-\frac{1}{2}&\frac{i\sqrt{3}}{2}\\
&\\
\frac{i\sqrt{3}}{2}&-\frac{1}{2}
\end{pmatrix}$$
And the matrix for a $\sigma_1$-braid is simply a diagonal matrix with phases on the diagonal.
Thus, an idea to create such ancillas is to start with the qubit $2211$ instead. We have seen when working on the qutrit that a full twist in the center has the effect of swapping $|0>$ and $|2>$. If instead we do a single braid in the center, we obtain the following matrix
$$\begin{array}{l}\qquad\qquad\qquad|0>\qquad|2>\\\\\begin{array}{cc}\begin{array}{l}|1>\\\\|3>\end{array}&\begin{pmatrix}
\frac{1}{\sqrt{2}}\,e^{i\frac{2\pi}{3}}&\frac{1}{\sqrt{2}}\\
&\\
\frac{1}{\sqrt{2}}\,e^{-i\frac{5\pi}{6}}&-\frac{i}{\sqrt{2}}
\end{pmatrix}\end{array}\end{array}$$
for the action

\begin{center}
\epsfig{file=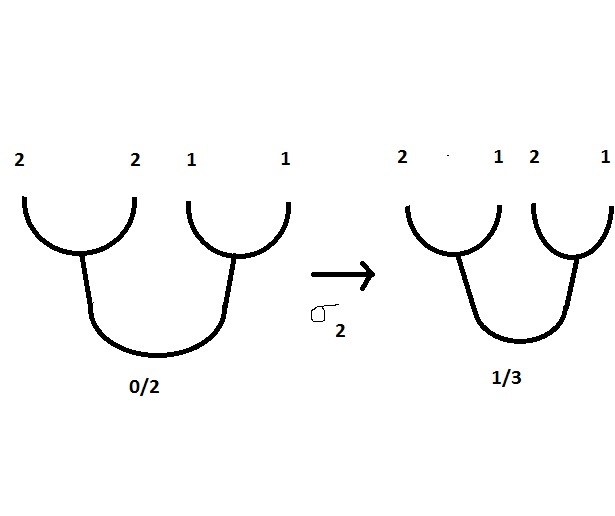, height=6cm}
\end{center}

Thus, by doing
\begin{center}\epsfig{file=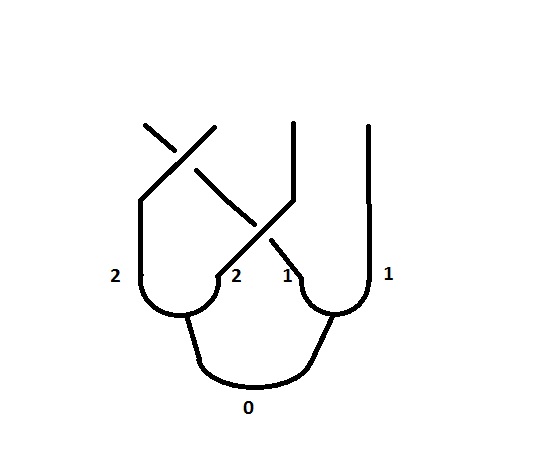, height=6cm}\end{center}
we obtain

\begin{center}
\epsfig{file=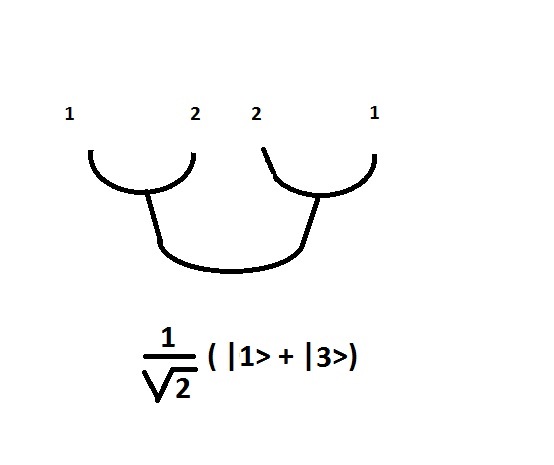, height=5cm}
\end{center}

Note that if you braid a $|2>$ instead, you can make the ancilla
$$\frac{1}{\sqrt{2}}(|1>-|3>)$$

\section{Discussion}

In the current paper, we give a physical interpretation of the actual original $D(18,1,1;2,1,1)$ as defined by generators in \cite{B},
while in \cite{LEV} we give a physical interpretation of an isomorphic copy of that group. \\%One of our main results is that the Blichfeld generator $F(18,1,1)$ can be realized physically by combining the $N$ gate and the FFO gate. \\
%The new set of "simpler looking" qutrit gates used here allows to uncover the structure of the group more easily than in our previous studies \cite{BL} and \cite{LEV}. Moreover, once we know the order of the group, we can easily derive the group is $D(18,1,1;2,1,1)$, without exhibiting complicated isomorphisms of groups like in \cite{LEV} or without the need of computer programs like GAP, also used in \cite{LEV} as a check-up. \\
It is disappointing but not surprising that we did not succeed to increase the number of qutrit gates by doing our protocol.
Enlarging such a number is not an easy problem. In fact, even using protocols with both braiding and interferometric measurement does not easily lead to finding additional gates which are not issued from braids we already have (cf Bauer's beautiful programming in \cite{F} to test such protocols by brute computer force). \\
%However, applying similar techniques yields a success for the qubit $1221$ when doing similar braids and several measurements.
%It is possible that the Blichfeld generators of $D(18,1,1;2,1,1)$ which we now know how to realize, could play a special role as ancillas leading to some interesting qubit gates obtained by braiding and forced measurement.
%The work of the present paper finds it origin in trying to solve the problem of leakage when doing measurements on some specific protocols which we are currently trying in a goal to find gates beyond braiding. Such gates would potentially enable universal quantum computation.
\begin{center} \textbf{Acknowledgements}\end{center}
%This work was done at and partially funded by Microsoft Research Station Q. The author thanks Station Q for hospitality and offering a quite exciting environment of work. She thanks Michael Freedman for numerous discussions.
This work is in the continuation of some previous works done at Microsoft Research Station Q.
The author is pleased to thank Michael Freedman for discussions.

\hfill\\
\noindent \textsc{Department of Mathematics, University of California, Santa Barbara, CA $93106$}\\\\
\textit{E-mail address:} \textbf{claire@math.ucsb.edu}

\end{document}